\documentclass[12pt,a4paper]{article}

\usepackage{amssymb}

\newcommand{\rset}{\mathbb{R}}
\newcommand{\nset}{\mathbb{N}}
\newcommand{\rmax}{\rset_{\max}}
\newcommand{\rmin}{\rset_{\min}}

\newcommand{\Mat}{\mathrm{Mat}}
\newcommand{\0}{\mathbf{0}}
\newcommand{\1}{\mathbf{1}}
\newcommand{\cle}{\preccurlyeq}
\newcommand{\clt}{\prec}
\newcommand{\cgt}{\succ}
\newcommand{\oodot}{\mathbin{\overline{\odot}}}
\newcommand{\ooplus}{\mathbin{\overline{\oplus}}}
\newcommand{\oostar}{\mathbin{\overline{\star}}}
\renewcommand{\a}{\mathbf a}

\newcommand{\A}{\mathbf A}
\newcommand{\B}{\mathbf B}
\newcommand{\E}{\mathbf E}
\newcommand{\I}{\mathbf I}
\newcommand{\lamint}{{[\underline\lambda,\overline\lambda]}}
\renewcommand{\t}{\mathbf t}
\newcommand{\V}{\mathbf V}
\newcommand{\x}{\mathbf x}
\newcommand{\X}{\mathbf X}

\newcommand{\y}{\mathbf y}
\newcommand{\z}{\mathbf z}
\newcommand{\lA}{\underline{\A}}

\newcommand{\lt}{\underline{\t}}
\newcommand{\lV}{\underline{\V}}
\newcommand{\lx}{\underline{\x}}
\newcommand{\ly}{\underline{\y}}
\newcommand{\lz}{\underline{\z}}

\newcommand{\uA}{\overline{\A}}

\newcommand{\ut}{\overline{\t}}
\newcommand{\uV}{\overline{\V}}
\newcommand{\ux}{\overline{\x}}
\newcommand{\uy}{\overline{\y}}
\newcommand{\uz}{\overline{\z}}

\newtheorem{thm}{Theorem}
\newtheorem{prop}{Proposition}

\newcounter{ExampleCount}[section]
\renewcommand{\theExampleCount}{\thesection.\arabic{ExampleCount}}
\newcommand{\example}[1]{%
\par%
\refstepcounter{ExampleCount}%
\textsc{Example \theExampleCount. #1{}}}

\newcounter{RemarkCount}[section]
\renewcommand{\theRemarkCount}{\thesection.\arabic{RemarkCount}}
\newcommand{\remark}{%
\par%
\refstepcounter{RemarkCount}%
\textsc{Remark \theRemarkCount. }}

\title{Idempotent Mathematics and Interval Analysis
\thanks{The work was supported by the joint
INTAS--RFBR grant No.\ 95--91 and the Erwin Schr\"odinger Institute
for Mathematical Physics (see Preprint ESI~632 at
\texttt{http://www.esi.ac.at}).}}

\author{G.~L.~Litvinov\thanks{Electronic mail: litvinov@islc.msk.su} \and
V.~P.~Maslov\thanks{Electronic mail: maslov@ipmnet.ru} \and
A.~N.~Sobolevski\u\i\thanks{Electronic mail: ansobol@idempan.phys.msu.su}}
\date{}

\begin{document}

\maketitle

\begin{abstract}
A brief introduction into Idempotent Mathematics and an idempotent version of
Interval Analysis are presented. Some applications are discussed.


\end{abstract}

\section{Introduction}%
\label{s:intro}

\quad\
Many problems in the optimization theory and other fields of
mathematics appear to be linear over semirings with idempotent addition (the
so-called \emph{idempotent superposition principle} \cite{Maslov}, which is a
natural analog of the well-known superposition principle in Quantum
Mechanics). The corresponding approach is developed systematically as
\emph{Idempotent Mathematics} or \emph{Idempotent Analysis}, a branch of
mathematics which has been growing vigorously last time (see, e.g.,
\cite{Maslov}--\cite{BCOQ}).

One of the most important examples of idempotent semirings is the set
$\rmax = \rset \cup \{-\infty\}$ equipped with operations $\oplus = \max$,
$\odot = +$ (see section~\ref{s:dequant}). In general, there exists a
correspondence between interesting, useful, and important constructions and
results concerning the field of real (or complex) numbers and similar
constructions dealing with various idempotent semirings. This
correspondence can be formulated in the spirit of the well known N.~Bohr's
correspondence principle in Quantum Mechanics; in fact, the two principles
are intimately connected (see \cite{LMCorrPrinc,LMShpiz,LMShpizDAN} and
sections~\ref{s:superpos}, \ref{s:matrices} and~\ref{s:correspondence}
below).  We discuss idempotent analogs of some basic ideas, constructions,
and results in traditional calculus and functional analysis; also, we show
that the correspondence principle is a powerful heuristic tool to apply
unexpected analogies and ideas borrowed from different areas of mathematics
(see, e.g., \cite{Maslov}--\cite{LMShpizDAN}).

The theory is well advanced and includes, in particular, new integration
theory, new linear algebra, spectral theory, and functional analysis. Its
applications include various optimization problems such as multicriteria
decision making, optimization on graphs, discrete optimization with a large
parameter (asymptotic problems), optimal design of computer systems and
computer media, optimal organization of parallel data processing, dynamic
programming, applications to differential equations, numerical analysis,
discrete event systems, computer science, discrete mathematics, mathematical
logic, etc.~(see, e.g., \cite{AdvSovMath}--\cite{LMRod} and references therein).

In section~\ref{s:dequant} we give a short heuristic introduction into
Idempotent Mathematics. Section~\ref{s:semirings} contains definitions of
basic concepts of idempotent arithmetic and several important examples. In
sections \ref{s:semimodules}--\ref{s:correspondence} we consider the notion
of linearity in Idempotent Analysis and indicate some of its applications
to idempotent linear algebra.

Due to imprecision of sources of input data in real-world problems, the
data usually come in a form of confidence intervals or other number sets
rather than exact quantities. Interval Analysis (see,
e.g.,~\cite{Kearfott}--\cite{Moore}) extends operations of traditional
calculus from numbers to number intervals, thus allowing to process such
imprecise data and control rounding errors in computations. To construct
the analog of calculus of intervals in the context of optimization theory
and Idempotent Analysis, we develop a set-valued extension of idempotent
arithmetic (see section~\ref{s:set-valued}).

The interval extension of an idempotent semiring is constructed in
sections~\ref{s:weak} and~\ref{s:stronger}. The idempotent interval
arithmetic appears to be remarkably simpler than its traditional analog.
For example, in the traditional interval arithmetic multiplication of
intervals is not distributive with respect to interval addition, while
idempotent interval arithmetic conserves distributivity.

A simple application of interval arithmetic to idempotent linear algebra is
discussed in section~\ref{s:appl}.  We stress that in the traditional
Interval Analysis the set of all square interval matrices of a given order
does not form a semigroup with respect to matrix multiplication: this
operation is not associative since distributivity is lost in traditional
interval arithmetic.  On the contrary, in the idempotent case associativity
is conserved.

Two properties that make the idempotent interval arithmetic so simple are
monotonicity of arithmetic operations and positivity of all elements of an
idempotent semiring. In general, idempotent interval analysis appears to be
best suited for treating the problems with order-preserving transformations
of imprecise data. We stress that this construction provides another example
of heuristic power of the idempotent correspondence principle.

\section{Dequantization and idempotent\newline correspondence principle}%
\label{s:dequant}

\quad\
Let $\rset$ be the field of real numbers and $\rset_+$ be the subset
of all non-negative numbers. Consider the following change of variable:  $$ u
\mapsto w = h \ln u, $$ where $u \in \rset_+ \setminus \{0\}$, $h > 0$; thus
$u = e^{w/h}$, $w \in \rset$. Denote by $\0$ the additional element $-\infty$
and by $S$ the extended real line $\rset \cup \{\0\}$. The above change of
variable has a natural extension $D_h$ to the whole $S$ by $D_h(0) = \0$;
also, we denote $D_h(1) = 0 = \1$.

Denote by $S_h$ the set $S$ equipped with the two operations $\oplus_h$
(generalized addition) and $\odot_h$ (generalized multiplication) such that
$D_h$ is a homomorphism of $\{\rset_+, +, \cdot\}$ to $\{S, \oplus_h,
\odot_h\}$. This means that $D_h(u_1 + u_2) = D_h(u_1) \oplus_h D_h(u_2)$ and
$D_h(u_1 \cdot u_2) = D_h(u_1) \odot_h D_h(u_2)$, i.e., $w_1 \odot_h w_2 =
w_1 + w_2$ and $w_1 \oplus_h w_2 = h \ln (e^{w_1/h} + e^{w_2/h})$.  It is
easy to prove that $w_1 \oplus_h w_2 \to \max\{w_1, w_2\}$ as $h \to 0$.

Denote by $\rmax$ the set $S = \rset \cup \{\0\}$ equipped with operations
$\oplus = \max$ and $\odot = +$, where $\0 = -\infty$, $\1 = 0$ as above.
Algebraic structures in $\rset_+$ and $S_h$ are isomorphic; therefore $\rmax$
is a result of a deformation of the structure in $\rset_+$.

We stress the obvious analogy with the quantization procedure, where $h$ is
the analog of the Planck constant. In these terms, $\rset_+$ (or $\rset$)
plays the part of a `quantum object' while $\rmax$ acts as a `classical' or
`semiclassical' object that arises as the result of a \emph{dequantization}
of this quantum object.

Likewise, denote by $\rmin$ the set $\rset \cup \{\0\}$ equipped with
operations $\oplus = \min$ and $\odot = +$, where $\0 = +\infty$ and $\1 =
0$. Clearly, the corresponding dequantization procedure is generated by the
change of variables $u \mapsto w = -h \ln u$.

Consider also the set $\rset \cup \{\0, \1\}$, where $\0 = -\infty$, $\1 =
+\infty$, together with the operations $\oplus = \max$ and $\odot=\min$.
Obviously, it can be obtained as a result of a `second dequantization' of
$\rset$ or $\rset_+$.

The algebras presented in this section are the most important examples of
idempotent semirings (for the general definition see section~\ref{s:semirings}), the central
algebraic structure of Idempotent Analysis. The basic object of the
traditional calculus is a \emph{function} defined on some set $X$ and taking
its values in the field $\rset$ (or $\mathbb{C}$); its idempotent analog is a
map $X \to S$, where $X$ is some set and $S = \rmin$, $\rmax$, or another
idempotent \emph{semiring}. Let us show that redefinition of basic
constructions of traditional calculus in terms of Idempotent Mathematics can
yield interesting and nontrivial results (see, e.g., \cite{Maslov}--\cite{Gunawardena} for
details and generalizations).

\example{Integration and measures.}%
\label{e:21}
To define an idempotent analog of the Riemann integral, consider a Riemann
sum for a function $\varphi(x)$, $x \in X = [a,b]$, and substitute semiring
operations $\oplus$ and $\odot$ for traditional addition and multiplication
of real numbers in its expression (for the sake of being definite, consider
the semiring $\rmax$):
$$
   \sum_{i = 1}^N \varphi(x_i) \cdot \Delta_i \quad\mapsto\quad
   \bigoplus_{i = 1}^N \varphi(x_i) \odot \Delta_i
   = \max_{i = 1, \ldots, N}\, (\varphi(x_i) + \Delta_i),
$$
where $a = x_0 < x_1 < \cdots < x_N = b$, $\Delta_i = x_i - x_{i - 1}$, $i =
1,\ldots,N$. As $\max_i \Delta_i \to 0$, the integral sum tends to
$$
   \int_X^\oplus \varphi(x)\, dx = \sup_{x \in X} \varphi(x)
$$
for any function $\varphi\colon X \to \rmax$ that is bounded. In general, the
set function
$$
   m_\varphi(B) = \sup_{x \in B} \varphi(x), \quad B \subset X,
$$
is called an \emph{$\rmax$-measure} on $X$; since $m_\varphi(\bigcup_\alpha
B_\alpha) = \sup_\alpha m_\varphi(B_\alpha)$, this measure is completely
additive. An idempotent integral with respect to this measure is defined as
$$
 \int_X^\oplus \psi(x)\, dm_\varphi
 = \int_X^\oplus \psi(x) \odot \varphi(x)\, dx
 = \sup_{x \in X}\, (\psi(x) + \varphi(x)).
$$

\example{Fourier--Legendre transform.}
Consider the topological
group $G = \rset^n$. The usual Fourier--Laplace transform is defined as
$$
   \varphi(x) \mapsto \widetilde\varphi(\xi)
   = \int_G e^{i\xi \cdot x} \varphi(x)\, dx,
$$
where $\exp(i\xi \cdot x)$ is a \emph{character} of the group $G$, i.e., a
solution of the following functional equation:
$$
   f(x + y) = f(x)f(y).
$$
The idempotent analog of this equation is
$$
   f(x + y) = f(x) \odot f(y) = f(x) + f(y).
$$
Hence `idempotent characters' of the group $G$ are linear functions of the
form $x \mapsto \xi \cdot x = \xi_1 x_1 + \cdots + \xi_n x_n$. Thus the
Fourier--Laplace transform turns into
$$
   \varphi(x) \mapsto \widetilde\varphi(\xi)
   = \int_G^\oplus \xi \cdot x \odot \varphi(x)\, dx
   = \sup_{x \in G}\, (\xi \cdot x + \varphi(x)).
$$
This equation differs from the well-known Legendre--Fenchel transform (see,
e.g.,~\cite{Rock}) in insignificant details.

These examples suggest the following formulation of the idempotent
correspondence principle \cite{LMCorrPrinc}:
\begin{quote}
\emph{There exists a heuristic correspondence between interesting, useful,
and important constructions and results over the field of real (or complex)
numbers and similar constructions and results over idempotent semirings in
the spirit of N. Bohr's correspondence principle in Quantum Mechanics.}
\end{quote}

So Idempotent Mathematics can be treated as a `classical shadow (or
counterpart)' of the traditional Mathematics over fields.

\section{Idempotent semirings: Basic definitions}%
\label{s:semirings}

\quad\
Consider a set $S$ equipped with two algebraic operations:
\emph{addition} $\oplus$ and \emph{multiplication} $\odot$. The triple
$\{S, \oplus, \odot\}$ is an \emph{idempotent semiring} if it satisfies
the following conditions (here and below, the symbol $\star$ denotes any of
the two operations $\oplus$, $\odot$):
\begin{itemize}
\item the addition $\oplus$ and the multiplication $\odot$ are
\emph{associative}: $ x \star (y \star z) = (x \star y) \star z$ for all $x,
y, z \in S$;
\item the addition $\oplus$ is \emph{commutative}: $x \oplus y = y \oplus x$
for all $x,y \in S$;
\item the addition $\oplus$ is \emph{idempotent}: $x \oplus x = x$ for all
$x\in S$;
\item the multiplication $\odot$ is \emph{distributive} with respect to the
addition $\oplus$: $x\odot(y\oplus z) = (x\odot y)\oplus(x\odot z)$ and
$(x\oplus y)\odot z = (x\odot z)\oplus(y\odot z)$ for all $x, y, z\in S$.
\end{itemize}
In the rest of this paper we shall sometimes drop the word `idempotent' when
the context is clear.

A \emph{unity} of an idempotent semiring $S$ is an element $\1 \in S$ such
that for all $x \in S$
$$
   \1 \odot x = x \odot \1 = x.
$$

A \emph{zero} of an idempotent semiring $S$ is an element $\0 \in S$ such
that $\0 \neq \1$ and for all $x \in S$
$$
 \0 \oplus x = x,\qquad \0 \odot x = x \odot \0 = \0.
$$

It is readily seen that if an idempotent semiring $S$ contains a unity
(a zero), then this unity (zero) is determined uniquely.

A semiring $S$ is said to be \emph{commutative} if $x \odot y = y \odot x$
for all $x, y \in S$. A commutative semiring is called a \emph{semifield}
if every nonzero element of this semiring is invertible. It is clear that
$\rmax$ and $\rmin$ are semifields.

Note that different versions of this axiomatics are used, see, e.g.,
\cite{AdvSovMath}--\cite{BCOQ} and some literature indicated in these books and papers.

The addition $\oplus$ defines on an idempotent semiring $S$ a~\emph{partial
order}:  $x \cle y$ iff $x \oplus y = y$. We use the notation $x \clt y$
if $x \cle y$ and $x \neq y$. If $S$ contains a zero $\0$, then $\0$ is its
least element with respect to the order $\cle$. The operations $\oplus$ and
$\odot$ are consistent with the order $\cle$ in the following sense: if
$x \cle y$, then $x \star z \cle y \star z$ and $z \star x \cle z \star y$
for all $x$, $y$, $z \in S$.

An idempotent semiring $S$ is said to be \emph{$a$-complete} if for any
subset $\{x_\alpha\} \subset S$, including $\varnothing$, there exists a sum
$\bigoplus\{x_\alpha\} = \bigoplus_\alpha x_\alpha$ such that
$(\bigoplus_\alpha x_\alpha) \odot y = \bigoplus_\alpha (x_\alpha \odot y)$
and $y \odot (\bigoplus_\alpha x_\alpha) = \bigoplus_\alpha (y \odot
x_\alpha)$ for any $y \in S$. An idempotent semiring $S$ containing a zero
$\0$ is said to be \emph{$b$-complete} if the conditions of
$a$-completeness are satisfied for any nonempty subset $\{x_\alpha\}
\subset S$ that is bounded from above. Any $b$-complete semiring either
is $a$-complete or becomes $a$-complete if the greatest element $\infty =
\sup S$ is added; see \cite{LMShpiz,LMShpizDAN} for details.

Note that $\bigoplus_\alpha x_\alpha = \sup\{x_\alpha\}$; in particular, an
$a$-complete idempotent semiring always contains the zero $\0 = \bigoplus
\varnothing$.

An idempotent semiring $S$ does not contain \emph{zero divisors} if $x \odot
y = \0$ implies that $x = \0$ or $y = \0$ for all $x$, $y \in S$.
An idempotent semiring $S$ is said to satisfy the~\emph{cancellation
condition} if for all $x$, $y$, $z \in S$ such that $x \neq \0$ it follows
from $x \odot y = x \odot z$ or $y \odot x = z \odot x$ that $y = z$. Any
idempotent semiring satisfying the cancellation condition does not contain
zero divisors. A commutative idempotent semiring $S$ is said to be
\emph{idempotent semifield} if every nonzero element of $S$ is invertible;
in this case the cancellation condition is fulfilled.

An idempotent semiring $S$ is said to be \emph{algebraically closed} if the
equation $x^n = y$, where $x^n = x\odot x\odot\cdots\odot x$ ($n$ times),
has a solution for all $y\in S$ and $n\in\nset$; an idempotent semiring $S$
with a unity $\1$ satisfies the~\emph{stabilization condition} if the
sequence $x^n \oplus y$ stabilizes whenever $x \cle \1$ and $y \neq
\0$~\cite{DSPreprint,DSAMS}. Note that in~\cite{DSAMS} the property of
algebraic closedness was incorrectly called `algebraic completeness'
due to translator's mistake.

The most important examples of idempotent semirings are those considered in
section~\ref{s:dequant}. We see that $\rmax$ is a $b$-complete algebraically closed
idempotent semifield satisfying stabilization condition. The idempotent
semiring $\rmin$ is isomorphic to $\rmax$. Note that both $\rmax$ and
$\rmin$ are linearly ordered with respect to the corresponding addition
operations; the order $\cle$ in $\rmax$ coincides with the usual linear
order $\leqslant$ and is opposite to the order $\cle$ in $\rmin$.

Consider the set $\widehat{\mathbb{R}}_{\max} = \rmax \cup \{\infty\}$ with
operations $\oplus$, $\odot$ extended by $\infty \oplus x = \infty$ for all
$x \in \rmax$, $\infty \odot x = \infty$ if $x \neq \0$ and $\infty \odot \0
= \0$. It is easily shown that this set is an $a$-complete idempotent
semiring and $\infty$ is its greatest element.

Let $\{S_1, S_2, \ldots\}$ be a collection if idempotent semirings. There
are several ways to construct a new idempotent semiring derived from the
semirings of this collection.

\example{}%
\label{e:31}
Suppose $S$ is an idempotent semiring and $X$ is an
arbitrary set. The set $\mathcal{M}(X;S)$ of all functions $X \to S$ is
an idempotent semiring with respect to the following operations:
$$
   (f \oplus g)(x) = f(x) \oplus g(x), \quad
   (f \odot g)(x) = f(x) \odot g(x), \quad x \in X.
$$
If $S$ contains a zero $\0$ and/or a unity $\1$, then the functions $o(x) =
\0$ for all $x \in X$, $e(x) = \1$ for all $x \in X$ are the zero and
the unity of the idempotent semiring $\mathcal{M}(X;S)$. It is also
possible to consider various subsemirings of $\mathcal{M}(X;S)$.

\example{}%
\label{e:32}
Let $S_i$ be idempotent semirings with operations
$\oplus_i$, $\odot_i$ and zeros $\0_i$, $i = 1, \ldots, n$. The set
$S = (S_1 \setminus \{\0_1\}) \times \cdots \times (S_n \setminus \{\0_n\})
\cup \0$ is an idempotent semiring with respect to the following operations:
$x \star y = (x_1, \ldots, x_n) \star (y_1, \ldots, y_n) = (x_1 \star_1 y_1,
\ldots, x_n \star_n y_n)$; the element $\0$ is the zero of this semiring.

Note that the set $\widetilde{S}=S_1\times\cdots\times S_n$ is also
an idempotent semiring with respect to the same operations; its zero is
the element $(\0_1,\ldots,\0_n)$.

Notice that even if primitive semirings in examples~\ref{e:31}
and~\ref{e:32} are linearly ordered sets with respect to the orders induced
by the correspondent addition operations, the derived semirings are only
partially ordered.

\example{}%
\label{e:33}
Let $S$ be an idempotent semiring and $\Mat_{mn}(S)$
be a set of all $S$-valued matrices. Define the sum $\oplus$ of matrices
$A = (a_{ij})$, $B = (b_{ij}) \in \Mat_{mn}(S)$ as $A \oplus B =
(a_{ij} \oplus b_{ij}) \in \Mat_{mn}(S)$, and let $\cle$ be the
corresponding order on $\Mat_{mn}(S)$. The \emph{product} of two matrices
$A \in \Mat_{lm}(S)$ and $B \in \Mat_{mn}(S)$ is a matrix $AB \in
\Mat_{ln}(S)$ such that $AB = (\bigoplus_{k = 1}^m a_{ik} \odot b_{kj})$.
The set $\Mat_{nn}(S)$ of square matrices is an idempotent semiring with
respect to these operations. If $\0$ is a zero of $S$, then the matrix $O =
(o_{ij})$, where $o_{ij} = \0$, is the zero of $\Mat_{nn}(S)$; if $\1$ is
a unity of $S$, then the matrix $E = (\delta_{ij})$, where $\delta_{ij} =
\1$ if $i = j$ and $\delta_{ij} = \0$ otherwise, is the unity of
$\Mat_{nn}(S)$.

Many additional examples can be found, e.g., in \cite{AdvSovMath}--\cite{BCOQ}.

\section{Idempotency and linearity}%
\label{s:semimodules}

\quad\
Now we discuss  an idempotent analog of a linear space. A set $V$  is
called a \emph{semimodule} over an idempotent semiring $S$ (or an
$S$-semimodule) if it is equipped with an idempotent commutative associative
addition operation $\oplus_V$ and a multiplication $\odot_V\colon S \times V
\to V$ satisfying the following conditions: for all $\lambda$, $\mu \in S$,
$v$, $w \in V$
\begin{itemize}
\item $(\lambda \odot \mu) \odot_V v = \lambda \odot_V (\mu \odot_V v)$;
\item $\lambda \odot_V (v \oplus_V w)
= (\lambda \odot_V v) \oplus_V (\lambda \odot_V w)$;
\item $(\lambda \oplus \mu) \odot_V v
= (\lambda \odot_V v) \oplus_V (\mu \odot_V v)$.
\end{itemize}
An $S$-semimodule $V$ is called a \emph{semimodule with zero} if $S$ is
a semiring with a zero $\0 \in S$ and there exists a zero element $\0_V \in
V$ such that for all $v \in V$, $\lambda \in S$
\begin{itemize}
\item $\0_V \oplus_V v = v$;
\item $\lambda \odot_V \0_V = \0 \odot_V v = \0_V$.
\end{itemize}

\example{Finitely generated free semimodule.}
The simplest
$S$-semimodule is the direct product $S^n = \{\, (a_1, \ldots, a_n) \mid a_j
\in S, j = 1, \ldots, n \,\}$. The set of all endomorphisms $S^n \to S^n$
coincides with the semiring $\Mat_{nn}(S)$ of all $S$-valued matrices (see
example~\ref{e:33}).

\example{Matrix semimodule.}
Take some $c \in S$, $A \in
\Mat_{mn}(S)$. The product $c \odot A$ is defined as the matrix $(c \odot
a_{ij}) \in \Mat_{mn}(S)$. Then the set of all $S$-valued matrices of a
given order $\Mat_{mn}(S)$ forms a semimodule under addition $\oplus$ and
multiplication by elements of $S$.

\example{Function spaces.}
An \emph{idempotent function space} $\mathcal{F}(X;S)$ is a subset of the
set of all maps $X \to S$ such that if $f(x)$, $g(x) \in \mathcal{F}(X;S)$
and $c \in S$, then $(f \oplus g)(x) = f(x) \oplus g(x) \in
\mathcal{F}(X;S)$ and $(c \odot f)(x) = c \odot f(x) \in \mathcal{F}(X;S)$;
thus an idempotent function space is another example of an $S$-semimodule.
If the semiring $S$ contains a zero element $\0$ and $\mathcal{F}(X;S)$
contains the zero constant function $o(x) = \0$, then the function space
$\mathcal{F}(X;S)$ has the structure of a semimodule with zero $o(x)$ over
the semiring $S$.

Recall that the idempotent addition defines a partial order in semiring $S$.
An important example of an idempotent functional space is the space
$\mathcal{B}(X;S)$ of all functions $X \to S$ bounded from above with respect
to the partial order $\cle$ in $S$. There are many interesting spaces of this
type including $\mathcal{C}(X;S)$ (a space of continuous functions defined on
a topological space $X$), analogs of the Sobolev spaces, etc.\ (see, e.g.,
\cite{AdvSovMath}--\cite{LMShpizDAN} for details).

According to the correspondence principle, many important concepts, ideas and
results can be converted from usual functional analysis to Idempotent
Analysis. For example, an idempotent scalar product in $\mathcal{B}(X;S)$ can
be defined by the formula
$$
   \langle\varphi,\psi\rangle = \int_X^\oplus \varphi(x) \odot \psi(x)\, dx,
$$
where the integral is defined as the $\sup$ operation (see example~\ref{e:21}).
Notice, however, that in the general case the ordering $\cle$ in $S$ is not
linear.

\example{Integral operators.}
It is natural to construct
idempotent analogs of \emph{integral operators} of the form
$$
   K:\, \varphi(y) \mapsto (K\varphi)(x)
   = \int_Y^\oplus K(x,y) \odot \varphi(y)\, dy,
$$
where $\varphi(y)$ is an element of a functional space $\mathcal{F}_1(Y;S)$,
$(K\varphi)(x)$ belongs to a space $\mathcal{F}_2(X;S)$ and $K(x,y)$ is a
function $X \times Y \to S$. Such operators are homomorphisms of the
corresponding functional semimodules. If $S = \rmax$, then this definition
turns into the formula
$$
   (K\varphi)(x) = \sup_{y \in Y}\, (K(x,y) + \varphi(y)).
$$
Formulas of this type are standard for optimization problems (see,
e.g.,~\cite{Bellman}).

\section{Idempotent superposition principle}%
\label{s:superpos}

\quad\
In Quantum Mechanics the superposition principle means that the
Schr\"odi\-n\-ger equation (which is basic for the theory) is linear.
Similarly in Idempotent Mathematics the (idempotent) superposition principle
means that some important and basic problems and equations (e.g.,
optimization problems, the Bellman equation and its versions and
generalizations, the Hamilton-Jacobi equation) nonlinear in the usual sense
can be  treated as linear over appropriate idempotent semirings, see
\cite{Maslov}--\cite{LMCorrPrinc}.

The linearity of the Hamilton-Jacobi equation over $\rset_{\min}$ (and
$\rset_{\max}$) can be deduced from the usual linearity (over $\mathbb{C}$)
of the corresponding Schr\"odinger equation by means of the dequantization
procedure described above (in section~\ref{s:dequant}). In this case the
parameter $h$ of this dequantization coincides with $i\hbar$ , where
$\hbar$ is the Planck constant; so in this case $\hbar$ must take imaginary
values (because $h>0$; see \cite{LMShpiz,LMShpizDAN} for details). Of
course, this is closely related to variational principles of mechanics; in
particular, the Feynman path integral representation of solution to
the Shr\"odinger equation corresponds to the Lax-Ole\u{\i}nik formula for
solution of the Hamilton-Jacobi equation (for the latter see, e.g.,
\cite{Lions}).

The situation is similar for the differential Bellman equation, see
\cite{Maslov,KolokolMaslov}.

B.A. Carr\'e \cite{Carre} used the idempotent linear algebra to show that
different optimization problems for finite graphs can be formulated in a
unified manner and reduced to solving Bellman equations, i.e., systems
of linear algebraic equations over idempotent semirings.

\textsc{Discrete Bellman equation.} It is well-known that discrete versions
of the Bellman equation can be treated as linear over appropriate idempotent
semirings. The following equation (the discrete stationary Bellman equation)
plays an important role in both discrete optimization theory and idempotent
matrix theory:
\begin{equation}
   X = AX \oplus B,
\label{Bellman}
\end{equation}
where $A \in \Mat_{nn}(S)$, $X, B \in \Mat_{nl}(S)$; matrices $A$, $B$ are
given and $X$ is unknown. Equation~(\ref{Bellman}) is a natural counterpart
of the usual linear system $AX = B$ in traditional linear algebra.

Note that if the closure matrix $A^*$ exists, then the matrix $X = A^* B$
satisfies~(\ref{Bellman}) because $A^* = AA^* \oplus E$. It can be
easily checked that this special solution is the minimal element of the set
of all solutions to~(\ref{Bellman}).

Actually, the theory of the discrete stationary Bellman equation
can be developed using the identity $A^* = AA^* \oplus E$ as an
additional axiom (the so-called \emph{closed semirings}; see, e.g.,
\cite{Lehmann}).

B.~A.~Carr\'e \cite{Carre} also generalized to the idempotent case
the principal algorithms of computational linear algebra and showed that
some of these coincide with algorithms independently developed for solution
of optimization problems; for example, Bellman's method of solving shortest
path problems corresponds to a version of Jacobi's method for solving
systems of linear equations, whereas Ford's algorithm corresponds to
a version of Gauss--Seidel's method.

We stress that these well-known results can be interpreted as a
manifestation of the idempotent superposition principle.

\section{Idempotent matrix theory: Some results}%
\label{s:matrices}

\textsc{Matrix algebra and graph theory.} Any square matrix $A = (a_{ij})
\in \Mat_{nn}(S)$ specifies a weighted directed graph with $n$ nodes such
that $a_{ij} \in S$ is the weight of the arc connecting $i$th node to $j$th
one. On the other hand, let $G$ be a directed graph with at most one arc
between any two nodes, every arc of which is characterized by its weight
that belongs to $S$. Then $G$ can be described by a matrix $A \in
\Mat_{nn}(S)$; in particular, nonexistent arcs are counted with weight $\0$,
and nonzero diagonal entries $a_{ii}$ of matrix $A$ correspond to loops.

For any $A \in \Mat_{nn}(S)$ define $A^0 = E$, $A^k = AA^{k - 1}$, $k \ge 1$.
Let $a^{(k)}_{ij}$ be $(i,j)$th element of the matrix $A^k$; then it is easy
to check that
$$
   a^{(k)}_{ij} =
   \bigoplus_{i_0 = i,\,
   1 \leqslant i_1, \ldots, i_{k - 1} \leqslant n,
   i_k = j} a_{i_0i_1} \odot \cdots \odot a_{i_{k - 1}i_k}.
$$
Consider a \emph{path} of length $k$ in the graph $G$, i.e. a sequence
of nodes $i_0,\dots,i_k$ connected by $k$ arcs with weights
$a_{i_0i_1},\dots,a_{i_{k - 1}i_k}$. The weight of the whole path is
defined to be the product $a_{i_0i_1} \odot \cdots \odot a_{i_{k - 1}i_k}$.
Thus $a^{(k)}_{ij}$ is the supremum of weights of all paths of length $k$
connecting node $i_0 = i$ to node $i_k = j$.

\textsc{Algebraic path problem.} This well-known problem is formulated as
follows: for each pair $(i,j)$ calculate supremum of weights of all paths
(of arbitrary length) connecting node $i$ to node $j$. If the semiring $S$
under investigation is $\rmin$ and arc weights are lengths in some metric,
then to solve this problem means to find all shortest paths. Is $S$ is
$\{\0, \1\}$ and the corresponding directed graph depicts some relation
$R$ in the set $\{1, \ldots, n\}$, then to solve this problem means to find
the transitive closure of $R$.

It is evident that in terms of matrix theory the algebraic path problem
is reduced to finding the matrix satisfying the formal expansion
$$
   A^* = E \oplus A \oplus A^2 \oplus \cdots = \bigoplus_{k = 0}^\infty A^k.
$$
The matrix $A^*$ is called the \emph{closure} of the matrix $A$.
If the idempotent semiring $S$ is not $a$-complete, this problem can be
nontrivial since we cannot simply take the infinite sum. Below we shall
discuss a sufficient condition for the existence of a closure, following
work of B.~A.~Carr\'e \cite{Carre}.

The matrix $A = (a_{ij}) \in \Mat_{nn}(S)$ is said to be \emph{definite}
(\emph{semi-definite}) if
$$
   a_{i_0i_1} \odot \cdots \odot a_{i_{k - 1}i_k} \clt \1
   \quad (a_{i_0i_1} \odot \cdots \odot a_{i_{k - 1}i_k} \cle \1)
$$
for any path $(i_1, \ldots, i_l)$ such that $i_0 = i_k$ (i.e., for any
closed path). Obviously, every definite matrix is semi-definite.
This definition is similar to that of B.~A.~Carr\'e~\cite{Carre} with
the only difference: Carr\'e considers an ordering that is opposite
to $\cle$.

\begin{thm}[Carr\'e]
Let $A$ be a semi-definite matrix. Then
$$
   \bigoplus_{l = 0}^k \A^l = \bigoplus_{l = 0}^{n - 1} \A^l
$$
for $k \geqslant n - 1$, so the closure matrix $A^* =
\bigoplus_{k = 0}^\infty A^k$ exists and is equal to
$\bigoplus_{k = 0}^{n - 1} A^k$.
\label{t:Carre}
\end{thm}
For the proof see \cite[Theorem~4.1]{Carre}. The basic idea of the proof
is evident: in the graph of a semi-definite matrix it is impossible to
construct a path of arbitrarily large weight since the weight of any closed
part of a path cannot be greater than $\1$. Thus there exists a universal
bound on path weights, which makes possible truncation of the infinite
series for the closure matrix.

\textsc{Spectral theory.} The spectral theory of matrices whose
elements lie in an idempotent semiring is similar to the well-known
Perron--Fro\-be\-ni\-us theory of nonnegative matrices (see,
e.g.,~\cite{KolokolMaslov,BCOQ,DSPreprint,DSAMS}).

Recall that the matrix $A = (a_{ij}) \in \Mat_{nn}(S)$
is said to be \emph{irreducible} in the sense of \cite{BCOQ} if for any
$1 \leqslant i,j \leqslant n$ either $a_{ij} \neq \0$ or there exist
$1 \leqslant i_1, \ldots, i_k \leqslant n$ such that $a_{ii_1} \odot
\cdots \odot a_{i_kj} \neq \0$. In~\cite{DSPreprint,DSAMS} matrices with
this property are called indecomposable.

We borrow the following important result from \cite{DSPreprint,DSAMS} (see
also \cite{BCOQ}):
\begin{thm}[Dudnikov, Samborski\u\i]
If a commutative idempotent semiring $S$ with a zero $\0$ and a unity $\1$
is algebraically closed and satisfies cancellation and stabilization
conditions, then for any matrix $A \in \Mat_{nn}(S)$ there exist a nonzero
`eigenvector' $V \in \Mat_{n1}(S)$ and an `eigenvalue' $\lambda \in S$ such
that $AV = \lambda \odot V$. If the matrix $A$ is irreducible, then the
`eigenvalue' $\lambda$ is determined uniquely.
\label{t:DSeigen}
\end{thm}
For the proof see \cite[Theorem 6.2]{DSAMS}.

An application of this result will be given in section~\ref{s:appl}.

Similar results hold for semimodules of bounded or continuous
functions~\cite{KolokolMaslov}.

\section{Correspondence principle for computations}%
\label{s:correspondence}

\quad\
Of course, the (idempotent) correspondence principle is valid for
algorithms as well as for their software and hardware implementations
\cite{LMCorrPrinc,LMRod}. Thus:
\begin{quote}
{\it If we have an important and interesting numerical algorithm, then there
is a good chance that its semiring analogs are important and interesting as
well.}
\end{quote}

In particular, according to the superposition principle, analogs of linear
algebra algorithms are especially important. Note that numerical algorithms
for standard infinite-dimensional linear problems over idempotent semirings
(i.e., for problems related to idempotent integration, integral operators and
transformations, the Hamilton-Jacobi and generalized Bellman equations) deal
with the corresponding finite-dimensional (or finite) ``linear
approximations''. Nonlinear algorithms often can be approximated by linear
ones. Thus idempotent linear algebra is the basis of the idempotent numerical
analysis.

Moreover, it is well-known that algorithms of linear algebra are convenient
for parallel computations; their idempotent analogs admit parallelization as
well. Thus we obtain a systematic way of applying parallel computation to
optimization problems.

Basic algorithms of linear algebra (such as inner product of two vectors,
matrix addition and multiplication, etc.) often do not depend on concrete
semirings, as well as on the nature of domains containing the elements of
vectors and matrices. Thus it seems reasonable to develop \emph{universal
algorithms} that can deal equally well with initial data of different domains
sharing the same basic structure~\cite{LMCorrPrinc,LMRod}; an example of such
universal Gauss--Jordan elimination algorithm is found in~\cite{Rote}.

Numerical algorithms are combinations of basic operations with `numbers',
which are elements of some numerical \emph{domains} (e.g., real numbers,
integers, etc.). But every computer uses some finite \emph{models} or finite
representations of these domains. Discrepancies between `ideal' numbers and
their `real' representations lead to calculation errors. This is another
reason to deal with universal algorithms that allow to choose a concrete
semiring and take into account the effects of its concrete finite
representation in a systematic way; see~\cite{LMCorrPrinc,LMRod} for details and
applications of the correspondence principle to hardware and software design.

\section{Set-valued idempotent arithmetic}%
\label{s:set-valued}

\quad\
Suppose $S$ is an idempotent semiring and $\mathcal{S}$ is a system of its
subsets. We shall denote the elements of $\mathcal{S}$ by $\x$, $\y$,~\dots
and define $\x \star \y = \{\, t \in S \mid t = x \star y, x \in \x, y \in
\y \,\}$.

We shall suppose that $\mathcal{S}$ satisfies the following two conditions:
\begin{itemize}
\item if $\x, \y \in \mathcal{S}$ and $\star$ is an algebraic operation in
$S$, then there exists $\z \in \mathcal{S}$ such that $\z \supset
\x \star \y$;
\item if $\{\x_\alpha\}$ is a subset of $\mathcal{S}$ such that
$\bigcap_\alpha \x_\alpha \neq \varnothing$, then there exists the infimum of
$\{\x_\alpha\}$ in $\mathcal{S}$ with respect to the ordering $\subset$,
i.e., the set $\y \in \mathcal{S}$ such that $\y \subset \bigcap_\alpha
\x_\alpha$ and $\z \subset \y$ for any $\z \in \mathcal{S}$ such that $\z
\subset \bigcap_\alpha \x_\alpha$.
\end{itemize}

Define algebraic operations $\ooplus$, $\oodot$ in
$\mathcal{S}$ as follows: if $\x, \y \in \mathcal{S}$, then $\x
\oostar \y$ is the infimum of the set of all elements $\z \in \mathcal{S}$
such that $\z \supset \x \star \y$.

\begin{prop}
The following assertions are true:
\begin{enumerate}
\item $\mathcal{S}$ is closed with respect to operations $\ooplus$,
$\oodot$.
\item The element $\x \oostar \y$ is optimal in the following sense: suppose
the exact values of input variables $x$ and $y$ lie in sets $\x$ and $\y$,
respectively; then the result of an algebraic operation $\x \oostar \y$
contains the quantity $x \star y$ and is the least subset of $S$ in
$\mathcal{S}$ with this property.
\item If the system $\mathcal{S}$ contains all one-element subsets of $S$,
then the semiring $\{S, \oplus, \odot\}$ is isomorphic to a subset of
the algebra $\{\mathcal{S}, \ooplus, \oodot\}$.
\end{enumerate}
\label{p:setvalued}
\end{prop}

The proof is straightforward.

In general, not much can be said about the algebra $\{\mathcal{S},
\ooplus, \oodot\}$, as the following example shows.

\example{}
Let $\mathcal{S} = 2^S$ and $\x \oostar \y = \x \star
\y$. In general, the set $\mathcal{S}$ with these `na\"{\i}ve' operations
$\ooplus$, $\oodot$ satisfies the above assumptions
but is not an idempotent semiring. Indeed, let $S$ be the semiring
$(\rmax\setminus\{\0\})^2 \cup\{(\0,\0)\}$ with operations $\oplus$, $\odot$
defined as in example~\ref{e:32}. Consider a set $\x = \{(0,1),(1,0)\} \in
\mathcal{S}$; we see that $\x \ooplus \x = \{(0,1),(1,0),(1,1)\}
\neq \x$ and if $\y = \{(1,0)\}$, $\z = \{(0,1)\}$, then $\x \oodot
(\y \ooplus \z) = \{(1,2),(2,1)\} \neq (\x \oodot \y)
\ooplus (\x \oodot z) = \{(1,1),(1,2),(2,1),(2,2)\}$. Thus
$\mathcal{S}$ with operations $\ooplus$, $\oodot$ does not
satisfy axioms of idempotency and distributivity.

It follows that $\mathcal{S}$ should satisfy some additional conditions in
order to have the structure of an idempotent semiring. In the next sections
we consider the case when $\mathcal{S}$ is a set of all closed intervals;
this case is of particular importance since it represents an idempotent
analog of the traditional Interval Analysis.

\section{Weak interval extension of an idempotent\newline semiring}%
\label{s:weak}

\quad\
Let $S$ be an idempotent semiring. A (closed) \emph{interval} in $S$
is a set of the form $\x = [\lx,{\ux}] = \{\, t \in S \mid \lx \cle t \cle
\ux \,\}$, where $\lx$, $\ux\in S$ ($\lx \cle \ux$) are said to be
\emph{lower} and \emph{upper bounds} of the interval $\x$, respectively.

Note that if $\x$ and $\y$ are intervals in $S$, then $\x \subset \y$ iff
$\ly \cle \lx \cle \ux \cle \uy$. In particular, $\x = \y$ iff $\lx = \ly$
and $\ux = \uy$.
\remark
Let $\x$, $\y$ be intervals in $S$. In general, the set $\x \star \y$ is not
an interval in $S$. Indeed, consider a set $S = \{\0,a,b,c,d\}$ and let
$\oplus$ be defined by the following order relation: $\0$ is the least
element, $d$ is the greatest element, and $a$, $b$, and $c$ are noncomparable
with each other. If $\odot$ is a zero multiplication, i.e., if $x \odot y =
\0$ for all $x$, $y \in S$, then $S$ is an idempotent semiring without unity.
Let $\x = [\0,a]$ and $\y = [\0,b]$; thus $\x \oplus \y = \{\0,a,b,d\}$. This
set is not an interval since it does not contain $c$ although $\0 \cle c
\cle d$.
\smallskip\par

Let $S$ be an idempotent semiring. A \emph{weak interval extension} $I(S)$ of
the semiring $S$ is the set of all intervals in $S$ equipped with
the operations $\ooplus$,~$\oodot$, where for any two
intervals $\x, \y \in I(S)$ $\x \oostar \y$ is defined as the smallest
interval containing the set $\x \star \y$.

\begin{prop}
$I(S)$ is closed with respect to the operations
$\ooplus$,~$\oodot$ and has the structure of an idempotent
semiring. Moreover, $\x \oostar \y = [\lx \star \ly, \ux \star \uy]$ for all
$\x, \y \in I(S)$.
\label{p:def_ostar}
\end{prop}
\textsc{Proof.} Take $t \in \x \star \y$ and let $x \in \x$, $y \in \y$ be
such that $t = x \star y$. By definition of interval, it follows that $\lx
\cle x \cle \ux$ and $\ly \cle y \cle \uy$.  Since the operation $\star$ is
consistent with the order $\cle$, we see that $\lx \star \ly \cle x \star y
\cle \ux \star \uy$. In particular, $\lx \star \ly \cle \ux \star \uy$, i.e.,
the interval $[\lx \star \ly, \ux \star \uy]$ is well defined. It follows
that $\x \star \y \subset [\lx \star \ly, \ux \star \uy]$.

Now let an interval $\z$ in $S$ be such that $\x \star \y \subset \z$. We
have $\lx \star \ly \in \x \star \y \subset \z$ and $\ux \star \uy \in \x
\star \y \subset \z$; hence $\lz \cle \lx \star \ly$ and $\ux \star \uy \cle
\uz$. Since $\lx \star \ly \cle \ux \star \uy$ by the above, it follows that
$[\lx \star \ly, \ux \star \uy] \subset \z$.

Thus the set $I(S)$ with the operations $\ooplus$,~$\oodot$
can be identified with a subset of an idempotent semiring $S \times S$ (see
example~\ref{e:32}). Since $\lx \star \ly \cle \ux \star \uy$ whenever $\lx \cle
\ux$ and $\ly \cle \uy$, $I(S)$ is closed with respect to
$\ooplus$,~$\oodot$; hence it is an idempotent semiring
(a subsemiring of $S\times S$). \hfill$\square$
\remark
Note that $I(S)$ satisfies the third condition of
proposition~\ref{p:setvalued} only if the semiring $S$ is $b$-complete;
nevertheless, the operations $\ooplus$,~$\oodot$ are
well defined in the general case.

\begin{prop}
If $S$ is an $a$-complete ($b$-complete) idempotent semiring, then $I(S)$
is an $a$-complete ($b$-complete) idempotent semiring.
\label{p:abcomplete}
\end{prop}
\textsc{Proof.} Let $S$ be an $a$-complete ($b$-complete) idempotent
semiring and $\{\x_\alpha\}$ be a nonempty subset of $I(S)$ (a nonempty
subset of $I(S)$ such that the set $\{\ux_\alpha\}$ is bounded in $S$). We
claim that the interval $\left[\bigoplus_\alpha \lx_\alpha, \bigoplus_\alpha
\ux_\alpha\right]$ contains the set $N = \{\, t \in S \mid (\forall\alpha)
(\exists x_\alpha \in \x_\alpha)\, t = \bigoplus_\alpha x_\alpha \,\}$ and
is contained in every other interval containing $N$; hence
$\overline{\bigoplus}_\alpha \x_\alpha = \left[\bigoplus_\alpha \lx_\alpha,
\bigoplus_\alpha \ux_\alpha \right]$. Indeed, in the case of an $a$-complete
semiring this statement is proved similarly to proposition~\ref{p:def_ostar}.
If $S$ is $b$-complete and the set $\{\ux_\alpha\}$ is bounded from above,
then the set $\{\lx_\alpha\}$ is also bounded from above, i.e., there exists
$y \in S$ such that $\lx_\alpha \cle \ux_\alpha \cle y$ for all $\alpha$.
Thus there exist $\bigoplus_\alpha\lx_\alpha \in S$ and $\bigoplus_\alpha
\ux_\alpha \in S$. Now the obvious adaptation of the above proof completes
the argument.

If $S$ is $a$-complete and $X \subset I(S)$ is empty, then by above
definitions $\overline{\bigoplus} X = [\0,\0]$ and for all $\y \in I(S)$ $\y
\oodot \left(\overline{\bigoplus} X \right) =
\left(\overline{\bigoplus} X \right) \oodot \y = [\0,\0]$. A direct
calculation shows that if $X = \{\x_\alpha\} \subset I(S)$ is nonempty, then
$$
   \y \oodot \left(\overline{\bigoplus}_\alpha \x_\alpha\right)
   = \overline{\bigoplus}_\alpha (\y \oodot \x_\alpha),
   \quad \left(\overline{\bigoplus}_\alpha \x_\alpha\right) \oodot \y
   = \overline{\bigoplus}_\alpha (\x_\alpha \oodot \y)
$$
for any $\y \in I(S)$. Thus $I(S)$ is $a$-complete ($b$-complete) if $S$ is
$a$-complete ($b$-complete). \hfill$\square$

The following two propositions are straightforward consequences of
proposition~\ref{p:def_ostar}:

\begin{prop}
If an idempotent semiring $S$ is commutative, so is $I(S)$.
\label{p:commute}
\end{prop}

\begin{prop}
If an idempotent semiring $S$ contains a zero $\0$ \emph{(\emph{unity}
$\1$)}, then the interval $[\0,\0]$ \emph{($[\1,\1]$)} is the zero
\emph{(\emph{unity})} of $I(S)$.
\label{p:zerounity}
\end{prop}

\begin{prop}
If an idempotent semiring $S$ has a zero $\0$ and does not contain zero
divisors, $I(S)$ has no zero divisors as well.
\label{p:zdiv}
\end{prop}
\textsc{Proof.} Let $\x, \y \in I(S)$ and $\x \neq [\0,\0]$, $\y \neq
[\0,\0]$. Because $\lx \cle \ux$, $\ly \cle \uy$, this means that $\ux \neq
\0$, $\uy \neq \0$. Thus if $\z = \x \oodot \y$, then $\uz = \ux
\odot \uy \neq \0$, since there are no zero divisors in $S$. It follows that
$\z \neq [\0,\0]$. \hfill$\square$

\begin{prop}
If $S$ is algebraically closed and for any $x, y \in S$, $n \in \nset$ the
equality $(x \oplus y)^n = x^n \oplus y^n$ holds, then $I(S)$ is
algebraically closed.
\label{p:aclosed}
\end{prop}
\textsc{Proof.} Suppose $\x^n = \x \oodot \x \oodot \cdots
\oodot \x = \y$. By proposition~\ref{p:def_ostar}, we see that $\lx^n
= \ly$ and $\ux^n = \uy$. Let $\lz \in S$ and $\uz \in S$ be the solutions of
these two equations. We claim that $\lz$ and $\uz$ can be chosen such that
$\lz \cle \uz$, i.e., the interval $[\lz,\uz]$ is well defined.

Take $\uz' = \lz \oplus \uz$; hence $\lz \cle \uz'$. Since in $S$ $\uz'^n =
(\lz \oplus \uz)^n = \lz^n \oplus \uz^n$, $\uz'^n = \ly \oplus \uy = \uy$. We
see that $[\lz,\uz']^n = [\ly,\uy] = \y$.  \hfill$\square$
\smallskip
\par\noindent
\remark%
\label{r:aclosed}
Note that the equality $(x \oplus y)^n = x^n \oplus y^n$
holds in any commutative idempotent semiring $S$ satisfying cancellation
condition (see, e.g., \cite{DSAMS}, assertion~2.1).

\section{A stronger notion of interval extension}%
\label{s:stronger}

We stress that in general a weak interval extension $I(S)$ of an idempotent
semiring $S$ that satisfies cancellation and stabilization conditions
does not inherit the latter two properties. Indeed, let $\x \oodot
\z = \y \oodot \z$, where $\z = [\0, z]$ and $z \neq \0$; then $\z$
is a nonzero element but this does not imply that $\x = \y$ since $\lx$ and
$\ly$ may not equal each other. Further, let $\y = [\0, y] \neq [\0, \0]$;
then the lower bound of $\x^n \overline \odot \y$ may not stabilize when
$n \to \infty$.

Therefore we define a stronger notion of \emph{interval extension} of an
idempotent semiring $S$ with a zero $\0$ to be the set $\I(S) = \{\, [x, y]
\mid x,y \in S,\, \0 \clt x \cle y, \,\} \cup \{[\0,\0]\}$ equipped with
operations $\ooplus$,~$\oodot$ defined as above.

Note that this object may not be well-defined if the semiring $S$ has zero
divisors. Indeed, let $\x = [x_1, x_2] \in \I(S)$ and $\y = [y_1, y_2] \in
\I(S)$ be such that $\0 \clt x_1 \clt x_2$, $\0 \clt y_1 \clt y_2$, $x_1
\odot y_1 = \0$, and $x_2 \odot y_2 \neq \0$; then $\x \oodot \y =
[\0, x_2 \odot y_2] \notin \I(S)$.

Throughout this section, we will suppose the interval extension $\I(S)$ of
an idempotent semiring $S$ to be closed with respect to the operations
$\ooplus$ and $\oodot$. To achieve this, it is sufficient to
require that the semiring $S$ contains no zero divisors.

\begin{thm}
The set $\I(S)$ is an idempotent semiring with respect to the operations
$\ooplus$ and $\oodot$ with the zero $[\0,\0]$ and does not
contain zero divisors. It inherits some special properties of the semiring
$S$:
\begin{enumerate}
\item If $S$ is $a$-complete ($b$-complete), then $\I(S)$ is $a$-complete
($b$-complete).
\item If $S$ is commutative, so is $\I(S)$.
\item If $\1$ is a unity of $S$, $[\1,\1]$ is the unity of $\I(S)$.
\item If $S$ is algebraically closed and for any $x, y \in S$, $n \in
\nset$ the equality $(x \oplus y)^n = x^n \oplus y^n$ holds, then $\I(S)$ is
algebraically closed.
\item If $S$ satisfies the cancellation condition, so does $\I(S)$.
\item If $S$ satisfies the stabilization condition, so does $\I(S)$.
\end{enumerate}
\label{t:Iproper}
\end{thm}

\textsc{Proof.} Using proposition~\ref{p:def_ostar}, it is easy to check that
$\I(S)$ is an idempotent semiring with respect to the operations
$\ooplus$,~$\oodot$. This semiring has the zero element
$[\0,\0]$ but contains no zero divisors by propositions~\ref{p:zerounity}
and~\ref{p:zdiv}. Propositions~\ref{p:abcomplete}--\ref{p:zerounity}
and~\ref{p:aclosed} imply the first four statements.

Suppose $S$ satisfies the cancellation condition, $\x$, $\y$, $\z \in \I(S)$,
and $\z$ is nonzero. If $\x \oodot \z = \y \oodot \z$, then
$\lx \odot \lz = \ly \odot \lz$ and $\ux \odot \uz = \uy \odot \uz$; since
$\z \neq [\0,\0]$ in $\I(S)$, $\lz \neq \0$ and $\uz \neq \0$, and it follows
from the assumptions that $\x = [\lx, \ux] = [\ly, \uy] = \y$. If $\z
\oodot \x = \z \oodot \y$, then $\x = \y$ similarly.

Suppose further that $S$ satisfies the stabilization condition; by
definition, $\ly \neq \0$ and $\uy \neq \0$ for any nonzero $\y \in \I(S)$.
Consider the sequence $\x^n \ooplus \y$; stabilization holds in $S$
for both bounds of the involved intervals and hence, by
proposition~\ref{p:def_ostar}, for the whole intervals as elements of
$\I(S)$. \hfill$\square$

Suppose $S$ is an idempotent semiring; then the map $\iota\colon S \to I(S)$
defined by $\iota(x) = [x,x]$ is an isomorphic imbedding of $S$ into its weak
interval extension $I(S)$. If the semiring $S$ has a zero $\0$ but no zero
divisors, then the map $\iota\colon S \to \I(S) \subset I(S)$ is an isomorphic
imbedding of $S$ into its interval extension. In the sequel, we will
identify the semiring $S$ with subsemirings $\iota(S) \subset I(S)$ or
$\iota(S) \subset \I(S) \subset I(S)$ and denote the operations in $I(S)$ or
$\I(S)$ by $\oplus$, $\odot$. If the semiring $S$ contains a unity $\1$,
then we denote the unit element $[\1,\1]$ of $I(S)$ or $\I(S)$ by $\1$;
similarly, we denote $[\0,\0]$ by $\0$.

\section{Cancellation and semifields}

We stress that in idempotent interval mathematics most of algebraic
properties of an idempotent semiring are conserved in its interval extension.
On the contrary, if $S$ is an idempotent semifield, then the set $\I(S)$ is
not a semifield but only a semiring satisfying the cancellation condition.

Recall that any commutative idempotent semiring $S$ with a zero $\0$ can be
isomorphically embedded into an idempotent semifield $\widetilde S$
provided that $S$ satisfies the cancellation condition (see,
e.g.,~\cite{DSPreprint}). If $\widetilde S$ coincides with its subsemifield
generated by $S$, then $\widetilde S$ is called a \emph{semifield of
fractions} corresponding to the semiring $S$. This semifield can be
constructed as the quotient $S \times (S \setminus \{\0\}) / \sim$, where
for any $(x, y)$,~$(z, t) \in S \times (S \setminus \{\0\})$
$$
   (x, y) \sim (z, t) \quad \Leftrightarrow \quad
   x \odot t = y \odot z,
$$
equipped with operations
$$
   (x, y) \oplus (z, t) = ((x \odot t) \oplus (y \odot z), y \odot t), \qquad
   (x, y) \odot (z, t) = (x \odot z, y \odot t).
$$
It is easy to see that these operations are defined in such a way that
pairs $(x, y)$ are treated as `fractions' with $x$ as `numerator' and $y$ as
(nonzero) `denominator'.  This semifield has the zero element $\{\, (\0, y)
\mid y \neq \0 \,\}$ and the unit element $\{\, (y, y) \mid y \neq \0
\,\}$; for every `fraction' $(x, y)$ representing a nonzero element of
$\widetilde S$ its inverse element is given by the fraction $(y, x)$.

In the context of the traditional Interval Analysis a similar extension of
the algebra of numerical intervals leads to the so-called \emph{Kaucher
interval arithmetic} \cite{Kaucher1,Kaucher2}. In addition to usual
intervals $[x, y]$, where $x \leqslant y$, it includes quasiintervals
$[x, y]$ with $y \leqslant x$, which arise as inverse elements for
the former.

Here we describe an idempotent version of Kaucher arithmetic. The following
statement shows that in this case the semifield of fractions of interval
extension $\I(S)$ corresponding to an idempotent semiring $S$ with
cancellation condition has very simple structure: it is isomorphic to the
idempotent semifield $(\widetilde S \setminus \{\0\})^2 \cup \{\0\}$.

\begin{prop}
Suppose $S$ is a commutative idempotent semiring with a zero $\0$, $S$
satisfies the cancellation condition, and $\widetilde S$ is its semifield of
fractions; then a semifield of fractions corresponding to the interval
extension $\I(S)$ is isomorphic to the semifield $(\widetilde S \setminus
\{\0\})^2 \cup \{\0\}$ (see example~\ref{e:32}).
\label{p:Grot}
\end{prop}
\textsc{Proof.} It follows from theorem~\ref{t:Iproper} that $\I(S)$ is a
commutative idempotent semiring with the zero element $\0 = [\0, \0]$ and
satisfies the cancellation condition. Thus $\I(S)$ can be isomorphically
embedded into its semifield of fractions.

Define the map $\phi\colon \I(S) \times (\I(S) \setminus \{\0\}) \to (\widetilde S
\setminus \{\0\})^2 \cup \{\0\}$ by the rule $\phi((\x,\y)) = (\lx \odot
\ly^{-1}, \ux \odot \uy^{-1})$. This map is surjective. Indeed, $(\0, \0) =
\phi((\0, \y))$ for any $\y \neq \0$; let us check that if $a, b \in \widetilde
S$, $a \neq \0$, $b \neq \0$, then there exist $\x, \y \in \I(S)$, $\y
\neq \0$, such that $(a, b) =  \phi((\x, \y))$. We see that there exist
$a_1, a_2, b_1, b_2 \in S$ such that $a = a_1 \odot a_2^{-1}$, $b =
b_1 \odot b_2^{-1}$. Define $\lx = a_1 \odot b_1 \odot b_2$, $\ux = (a_1
\odot b_1 \odot b_2) \oplus (a_2 \odot b_1^2)$, $\ly = a_2 \odot b_1 \odot
b_2$, $\uy = (a_1 \odot b_2^2) \oplus (a_2 \odot b_1 \odot b_2)$; thus $\lx
\cle \ux$, $\ly \cle \uy$ and $\lx \odot \ly^{-1} = a_1 \odot a_2^{-1} =
a$, $\ux \odot \uy^{-1} = b_1 \odot b_2^{-1} = b$.

Since $x \odot y^{-1} = z \odot t^{-1}$ iff $x \odot t = y \odot z$ for any
$x, y, z, t \in \widetilde S$, we see that $\phi((\x, \y)) = \phi((\z,
\t))$ iff $(\x, \y)$ and $(\z, \t)$ define the same element of the
semifield of fractions corresponding to $\I(S)$. Also,
$$
\begin{array}{rl}
   \phi(((\x \odot \t) \oplus (\y \odot \z), \y \odot \t))
   &\!\!\!\!= ((\lx \odot \ly^{-1})) \oplus (\lz \odot \lt^{-1}),
   (\ux \odot \uy^{-1})) \oplus (\uz \odot \ut^{-1})) \\
   &\!\!\!\!= \phi((\x, \y)) \oplus \phi((\z, \t)),\\
   \phi((\x \odot \z, \y \odot \t))
   &\!\!\!\!= ((\lx \odot \ly^{-1}) \odot (\lz \odot \lt^{-1}),
   (\ux \odot \uy^{-1}) \odot (\uz \odot \ut^{-1})) \\
   &\!\!\!\!= \phi((\x, \y)) \odot \phi((\z, \t)).\\
\end{array}
$$
Thus the semifield of fractions corresponding to $\I(S)$ is isomorphic to
the idempotent semifield $(\widetilde S \setminus \{\0\})^2 \cup \{\0\}$.
\hfill$\square$

The commutativity condition in this proposition is a natural one. Indeed,
it can be proved that if each nonzero element of a $b$-complete idempotent
semigroup $S$ has a multiplicative inverse, then $S$ is commutative (see,
e.g., \cite{LMShpiz}).

\section{Application to linear algebra}%
\label{s:appl}

\quad\
Suppose $S$ is an idempotent semiring with a zero $\0$ and a unity $\1$ and
$I(S)$ is its weak interval extension; then $\Mat_{nn}(I(S))$ is an
idempotent semiring. If the interval extension $\I(S)$ of the semiring $S$
is well defined, then the same is true for $\Mat_{nn}(\I(S))$. We shall
denote the (common) unit element of these semirings by $\E$.

If $\A = (\a_{ij}) \in \Mat_{mn}(I(S))$ ($\A = (\a_{ij}) \in
\Mat_{mn}(\I(S))$) is a (not necessarily square) interval matrix, then the
matrices $\lA = (\underline{\a_{ij}})$ and $\uA = (\overline{\a_{ij}})$ are
called \emph{lower} and \emph{upper matrices} of the \emph{interval matrix}
$\A$.

\begin{prop}
Let $S$ be an idempotent semiring with a zero $\0$ and a unity $\1$. The
mapping $\A \in \Mat_{nn}(I(S)) \mapsto [\lA, \uA] \in I(\Mat_{nn}(S))$ is
an isomorphism of idempotent semirings $\Mat_{nn}(I(S))$ and
$I(\Mat_{nn}(S))$. If the semiring $S$ has an interval extension $\I(S)$,
then this assertion remains true if $I(S)$ is substituted by
$\mathbf{I}(S)$.
\label{p:IMatMatI}
\end{prop}

Here intervals $[\lA, \uA]$ in $I(\Mat_{nn}(S))$ or $\I(\Mat_{nn}(S))$ are
defined with respect to the partial ordering $\cle$ in $\Mat_{nn}(S)$ (see
example~\ref{e:33}). The proof follows easily from
proposition~\ref{p:def_ostar}; indeed, this proposition implies that
addition and multiplication of interval matrices are reduced to separate
addition and multiplication of their lower and upper matrices.

The following proposition is an immediate consequence of
theorem~\ref{t:DSeigen}:
\begin{prop}
If a commutative idempotent semiring $S$ with a zero $\0$ and a unity $\1$
is algebraically closed and satisfies cancellation and stabilization
conditions, then for any matrix $\A \in \Mat_{nn}(\I(S))$ there exist a
nonzero `eigenvector' $\V \in \Mat_{n1}(\I(S))$ and an `eigenvalue'
$\lamint \in \I(S)$ such that $\A\V = \lamint \odot \V$. If the matrix $\A$
is irreducible, then the `eigenvalue' $\lamint$ is determined uniquely.
\label{p:DSeigen}
\end{prop}

It follows from proposition~\ref{p:def_ostar} that $\lA\,\lV =
\underline{\lambda} \odot \lV$ and $\uA\,\uV = \overline\lambda \odot \uV$.

Consider the following interval discrete stationary Bellman equation (see
also sections \ref{s:superpos},~ref{s:matrices}):
$$
   \X = \A\X \oplus \B,
$$
where $\A \in \Mat_{nn}(\I(S))$, $\B, \X \in \Mat_{ns}(\I(S))$.  Consider
the following iterative process:
\begin{equation}
   \X_{k + 1} = \A\X_k \oplus \B
   = \A^{k + 1}\X_0 \oplus \left(\bigoplus_{l = 0}^k \A^l\right) \B,
\label{iterate}
\end{equation}
where $\X_k \in \Mat_{ns}(\I(S))$, $k = 0, 1, \ldots$

The following proposition is due to B.~A.~Carr\'e \cite[Theorem~6.1]{Carre}
up to some terminology:
\begin{prop}
If matrix $A \in Mat_{nn}(S)$ is semi-definite, then the iterative process
$X_{k + 1} = AX_k \oplus B$ stabilizes to the solution $X = A^* B$ of the
equation $X = AX \oplus B$ after at most $n$ iterations for any initial
approximation $X_0 \in \Mat_{n1}(S)$ such that $X_0 \cle A^* B$.
\label{p:Carre}
\end{prop}

Suppose an idempotent semiring $S$ satisfies the assumptions of
proposition~\ref{p:DSeigen}. Let $[\underline{\lambda}_1,
\overline{\lambda}_1], \ldots, \allowbreak [\underline{\lambda}_k,
\overline{\lambda}_k]$, $1 \leqslant k \leqslant n$, be the eigenvalues of
the matrix $\A \in \Mat_{nn}(\I(S))$. Denote
$\sup\{\overline{\lambda}_1,\ldots,\overline{\lambda}_k\} =
\bigoplus_{l = 1}^k \overline{\lambda}_l$ by $\rho(\A)$. It is possible to
give a simple spectral criterion of convergence of the iterative
process~(\ref{iterate}):

\begin{thm}
Let $S$ be a commutative semiring satisfying conditions of
proposition~\ref{p:DSeigen} and matrix $\A \in \Mat_{nn}(\I(S))$ be such that
$\rho(\A) \cle \1$. Then the iterative process $\X_{k + 1} = \A\X_k
\oplus \B$, $k \geqslant 0$, stabilizes to the minimal solution $\X =
\A^* \B$ of equation $\X = \A\X \oplus \B$ after at most $n$ iterations for
any $\X_0 \in \Mat_{n1}(\I(S))$ such that $\X_0 \cle \X$.
\label{p:spectral}
\end{thm}
\textsc{Proof.} It follows from proposition~\ref{p:def_ostar} that it is
sufficient to prove that sequences of lower and upper matrices of
$\{\X_k\}$ converge separately. To this end, we shall show that the matrices
$\lA$ and $\uA$ are semi-definite; then the result will follow from
proposition~\ref{p:Carre}.

Since $\underline{\a_{ij}} \cle \overline{\a_{ij}}$ for all $i,j$, we need
only to prove that $\uA$ is semi-definite. First we shall prove this if
$\uA$ is irreducible. Using the expression for a unique eigenvalue of an
irreducible matrix $\uA$ in terms of cycle invariants~\cite{DSAMS}
$$
   \overline\lambda^{\phi(n)} =
   \bigoplus_{l = 1, \ldots, n \atop (i_1, \ldots, i_l)}
   \left[ \overline{\a_{i_1i_2}} \odot \cdots \odot
   \overline{\a_{i_li_1}} \right]^{\phi(n)/l},
$$
where $\phi(n)$ is the least common multiple of the numbers $1, \ldots, n$,
we see that for any cycle its cycle invariant $P = \overline{\a_{i_1i_2}}
\odot \cdots \odot \overline{\a_{i_li_1}}$ satisfies $P \cle \1$ if
$\overline\lambda \cle \1$ (indeed, if $P \oplus \1 \cgt \1$, then, by
remark~\ref{r:aclosed}, $(\1 \oplus P)^{\phi(n)/l} =
\1 \oplus P^{\phi(n)/l} \cgt \1$, so $\1 \oplus \overline\lambda^{\phi(n)}
\cgt \1$---a contradiction). Thus $\uA$ is a semi-definite matrix.

If $\uA$ is reducible, there exists a permutation matrix $Q$ such that $\uA =
Q B Q^{-1}$, where $B = (b_{ij})_{i,j = 1,\ldots,n}$ has an upper
block triangular form:
$$
   B = \left(\matrix{B_1    &  *     & \cdots & *      \cr
                     \0     & B_2    & \cdots & *      \cr
                     \cdots & \cdots & \cdots & \cdots \cr
                     \0     & \0     & \cdots & B_k    \cr}\right),\qquad
   1 < k \leqslant n,
$$
and all square matrices $B_1, \ldots, B_k$ are either zero or irreducible.
We claim that every eigenvalue $\overline\lambda$ of $\uA$ is an eigenvalue
of some $B_l$, $l = 1, \ldots, k$. Indeed, let $V$ be an eigenvector of
$\uA$ with an eigenvalue $\overline\lambda$; denote $i$th element of the
vector $V$ by $v_i$. Consider a decomposition $\{1, \ldots, n\} = N_1 \cup
\cdots \cup N_k$, where $N_i \cap N_j = \varnothing$ if $i \neq j$ and $B_l
= (b_{ij})_{i,j \in N_l}$, $l = 1, \ldots, k$; let $l_0 = \max \{\, l
\mid \exists i \in N_l\colon v_i \neq \0 \,\}$. We see that $\lambda$ is a
unique eigenvalue of the irreducible matrix $B_{l_0}$ corresponding to the
eigenvector $(v_i)_{i \in N_{l_0}}$. The condition $\rho(\A) \cle \1$
implies that $B_1, \ldots, B_k$ are semi-definite. Since $P = \0 \clt \1$
for any cycle containing indices $i \in N_l$, $ j \in N_s$, $l \neq s$, we
conclude that $\uA$ is a semi-definite matrix.  \hfill$\square$
%
\smallskip
\par\noindent
\remark
Compare this simple proposition with the well-known spectral
convergence criterion of the iterative process in traditional Interval
Analysis (\cite[theorem~12.1]{AlefHerz}), which in our notation has the following
form:
\smallskip
\par\noindent\emph{The iterative process $\X_{k + 1} = \A\X_k + \B$, $k \ge
0$, converges to a unique solution $\X$ of the equation $\X = \A\X + \B$
for any $\X_0 \in \Mat_{n1}(I(\mathbb{C}))$ if and only if $\rho(|\A|) <
1$.}

International Sophus Lie Centre,

Nagornaya, 27--4--72, Moscow 113186 Russia

\end{document}